\documentclass[letterpaper, 10 pt, conference]{ieeeconf}

\usepackage{times}

\usepackage{multicol}
\usepackage[bookmarks=true]{hyperref}
\usepackage{graphicx}
\usepackage{amsmath}
\usepackage{amssymb}
\usepackage{multirow}

\usepackage{xargs}

\IEEEoverridecommandlockouts                              

\title{\Large \bf
Evaluating direct transcription and nonlinear optimization methods for robot
motion planning
}

\author{Diego Pardo, Lukas M\"{o}ller, Michael Neunert, Alexander W. Winkler and
Jonas Buchli%
\thanks{
This research has been funded through a Swiss National Science Foundation
Professorship award to Jonas Buchli and by the Swiss National Centre of
Competence in Research NCCR - Robotics.
Diego Pardo \{{\tt\small depardo@ethz.ch}\},
Lukas M\"{o}ller \{{\tt\small moellelu@student.ethz.ch}\},
Michael Neunert \{{\tt\small neunertm@ethz.ch}\},
Alexander Winkler \{{\tt\small winklera@ethz.ch}\},
Jonas Buchli \{{\tt\small buchlij@ethz.ch}\}
are with the Agile \& Dexterous Robotics Lab at the 
Institute of Robotics and Intelligent Systems, ETH Z{\"u}rich, Switzerland. A
preprint of this paper is available at Arxiv:1504.05803. This manuscript
includes additional results.}
}

\begin{document}

\maketitle
\thispagestyle{empty}
\pagestyle{empty}

\begin{abstract}
This paper studies existing direct transcription methods for trajectory
optimization applied to robot motion planning.
%
%
%
%
%
%
There are diverse alternatives for the implementation of direct transcription. 
In this study we analyze the effects of such alternatives when solving a
robotics problem.
Different parameters such as integration scheme, number of
discretization nodes, initialization strategies and complexity of the problem
are evaluated.
We measure the performance of the methods in terms of computational time,
accuracy and quality of the solution. 
Additionally, we compare two optimization methodologies frequently used to
solve the transcribed problem, namely Sequential
Quadratic Programming (SQP) and Interior Point Method (IPM). 
As a benchmark, we solve different motion tasks on an underactuated and
non-minimal-phase ball-balancing robot with a 10 dimensional state space and 3
dimensional input space.
Additionally, we validate the results on a simulated 3D
quadrotor.
Finally, as a verification of using direct transcription methods for trajectory
optimization on real robots, we present hardware experiments on a motion task
including path constraints and actuation limits.
%

\end{abstract}

\IEEEpeerreviewmaketitle

\section{Introduction}

Numerical methods for trajectory optimization (TO) \cite{Betts1998} have
received considerable attention from the robotics community in recent years.
These methods are especially convenient for motion planning and control problems
on high dimensional systems with complex dynamics \cite{Koenemann2015,Dai2014}.
In such systems, using common approaches separating kinematic planning and control
struggles in finding plausible solutions and tends to produce inefficient and
artificial motions.
In contrast to designed or motion captured references, TO methods
generate dynamically feasible motions maximizing a performance criterion while
satisfying a set of constraints.
This is an appealing approach when dealing with nonlinear and unstable systems
like legged or balancing robots, where an optimal  solution cannot be obtained
analytically from the necessary conditions of optimality.

There are different types of methods to numerically solve the TO problem
(direct, indirect, single shooting, multiple shooting). 
The classification of these methods is out of the scope of this paper \cite{Betts1998}.
%
%
Here we are interested in studying the performance of the family of methods known
as \textsl{direct transcription} when applied to systems with high
nonlinearities and unstable dynamics.
Recent results in the area of whole-body dynamic motion planning have
demonstrated the potential of these methods \cite{Posa2014,Schultz2010,Dai2014},
where very dynamic and complex motions have been obtained.
Interestingly, to the best of our knowledge there are no reports of hardware
experiments of direct transcription methods on floating-base and naturally
unstable robots. 
%

In direct transcription the state and control trajectories are discretized over
time,
and an augmented version of the original problem is stated over the values of
the trajectories at the discrete points or nodes
\cite{Hargraves1987,vonStryk1992}.
Therefore, constraints enforcing the system dynamics between nodes need to be
included.
The resulting problem can be solved by a nonlinear programming (NLP) solver.
There are many approaches to direct transcription, mainly differing in the
method used to enforce the dynamic constraints.

During robot motion planning, the procedure consists of
folding the problem into
a continuous TO framework and subsequently find optimal trajectories 
using direct transcription and an NLP solver.
The main challenge of TO for robotics applications is to adequately integrate
all components of the problem into this formulation
(e.g. contact forces \cite{Dai2014}, contact dynamics \cite{Schultz2010}).
%
%
However, the transcription and solver stages are often treated as black box
processes with few degrees of freedom for the robotics researcher.
Specialized software is available to apply direct transcription on general
optimal control problems (e.g., PSOPT \cite{Becerra2010}), and
a number of NLP solvers exist (e.g. SNOPT \cite{Gill2005}, IPOPT
\cite{Wachter2006}), based on mature methods in mathematical optimization.

One of the main disadvantages of direct transcription is the augmented size of
the subsequent NLP, making its solution computationally expensive. 
This is especially true for robotics problems with high dimensional systems
with contact forces or obstacle constraints.
%
%
The accuracy of the solution is also a major concern, as the dynamics of the
robot are approximated in between the nodes and the resulting plan might not be
feasible on the real robot. Given such difficulties, it is worth to evaluate
the impact of choosing different types of dynamic constraints, number of
discretization nodes and NLP solvers on the accuracy and quality of the 
solution.
Moreover, it is also important to measure the performance of the different
implementations for tasks with diverse levels of complexity. 
%


%
%
The main contribution of this paper is the quantification of the performance of
various direct transcription methods in an experimental study. It is also shown
that these methods are sufficiently fast so online motion
planning on robotic hardware comes in reach.
This work also investigates the performance of two types of methods,
namely Sequential Quadratic Programming (SQP) and Interior Point Methods (IPM),
for solving the corresponding NLP problem. 
This paper does not present a survey or comparison of methods for solving TO
problems \cite{Conway2012}, nor does it evaluate available optimization software
packages or compare NLP solvers \cite{Betts2002}.
%
Instead, we use an unstable ball-balancing robot \cite{Fankhauser2010} (10
states, 3 inputs)
to compare the performance of these techniques in terms of computational time,
accuracy and quality of the solution.  
%
%
Influences of the complexity of the task and NLP initialization are also
evaluated.
Moreover, the results of the evaluation are validated using a different robot
and task (3D quadrotor, 12 states, 4 controls), 
indicating how results scale to other systems.


%

%
This paper is organized as follows. Section \ref{sec:direct_transcription}
presents a theoretical review of direct transcription methods and NLP solvers.
In Section \ref{sec:problem_set}, the methodology designed to evaluate and
compare the methods is described. Information about the model of the robot used
in this study, together with a description of the cost function and feedback
controller used during the benchmark tasks are provided in Section
\ref{sec:robot_model_and_control_scheme}. Results are shown in Section
\ref{sec:results} and analyzed and discussed in Section \ref{sec:discussion}.
Section \ref{sec:conclusions} presents the conclusions.

\section{Background}
\label{sec:direct_transcription}

In this section we describe direct transcription as well as
the key concepts supporting SQP and IPM solvers.
\subsection{Direct transcription for trajectory optimization}

In general, the system dynamics of a nonlinear robot can be modeled by a set of
differential equations,
\begin{equation}
\dot{x}(t) = f(x(t),u(t)),
\label{eq:systemdynamics}
\end{equation}
where $x \in \mathbb{R}^n$ represents the system states and $u \in \mathbb{R}^m$
the vector of control actions. The transition function $f(\cdot)$ defines the
system evolution in time. A single phase TO problem consists in finding a
finite-time input trajectory $u(t), \forall t \in [0,T]$,
such that a given criteria is minimized,
\begin{equation}
J = \Psi(x(T)) + \int_{0}^{T} \psi(x(t),u(t),t)~dt
\label{eq:cost_function}
\end{equation}
where $\psi(\cdot)$ and $\Psi(\cdot)$ are the intermediate and final cost
functions respectively. The optimization may be subject to a set of boundary and
path constraints,
%
\begin{eqnarray}
	\label{eq:constraint_3}
	\phi_{p,min}  \leq  \phi_p(x(t),u(t),t)  \leq  \phi_{p,max}
\end{eqnarray}
and bounds on the state and control variables
\begin{eqnarray}
x_{min} \leq x(t) \leq x_{max} \ \ , \ \
u_{min} \leq u(t) \leq u_{max}.
\label{eq:constraint_5}
\end{eqnarray}

Direct transcription translates this continuous formulation into an optimization
problem with a finite number of variables.
The set of decision variables, $y\in \mathbb{R}^p$, includes the discrete values
of the state and control trajectories at certain points or \textsl{nodes}.
Therefore,  $y=\{x_k,u_k\}$, for $k = 1,...,N$.
Moreover, the set of decision variables can be augmented with additional
parameters to be optimized. 
For instance, in the results presented in this paper the time between nodes
$\Delta T = t_{k+1}-t_k$ has been included as a decision variable.

The resulting NLP is then formulated as follows,
\begin{equation}
	\begin{array}{rlclcl}
		\displaystyle \min_{y} & f_0(y) \\
		\textrm{s.t.} & & \zeta(y) & =&  0& \\ 
		&  g_{min}  &\leq &g(y)& \leq & g_{max} \\
		& y_{min}  &\leq &y& \leq  &y_{max}, \\		
	\end{array}
\label{eq:NLP}
\end{equation}
where $f_0(\cdot)$ is a scalar objective function which in our implementation is
given by a quadrature formula approximating (\ref{eq:cost_function}), whereas 
the boundary and path constraints in (\ref{eq:constraint_3}) 
are gathered in $g(\cdot)$. The NLP also includes bounds on the decision
variables.
Additionally, a vector of dynamic constraints or \textsl{defects} $\zeta(\cdot)
\in \mathbb{R}^{(N-1)n}$ is added to verify the system dynamics at each
interval.

As a consequence of the discretization, the resulting NLP is considerably large.
This fact is partially compensated by the sparsity of the resulting problem, as
such structural property is handled particularly well by large-scale NLP solvers
\cite{Gill2005,Wachter2006}.

Direct transcription methods mainly differ in the way they formulate
the dynamic constraints in (\ref{eq:NLP}). 
The simplest dynamic constraint is given by Euler's integration rule, 
\begin{equation}
\zeta_k = x_{k+1} - x_k - f(x_k,u_k)\Delta T = 0.
\label{eq:euler}
\end{equation}
There are other approaches using different implicit integration rules.
For example, using a trapezoidal \cite{Betts1998} interpolation for the
dynamics, the defects are then given by
\begin{equation}
x_{k+1} - x_k - \frac{\Delta T}{2}\left[ f(t_k) + f(t_{k+1}) \right] = 
0,
\label{eq:trapezoidal}
\end{equation}
where the notation $f(t_i) = f(x_i,u_i)$ has been adopted for simplicity.
In the original formulation of the scheme known as \textsl{direct collocation}
\cite{vonStryk1992} \cite{Hargraves1987} piecewise cubic Hermite functions are
used to interpolate the states between nodes.
In this method the difference between the system dynamics and the derivative of
the Hermite function at the middle of the interval, $(t = t_c)$, is used as
dynamic constraint. 
%
%
%
This approach is equivalent to a Simpson's integration rule,
\begin{equation}
x_{k+1} - x_{k} - \frac{\Delta T}{6} \left[ f(t_k) + 4f(t_c) +
f(t_{k+1})\right] = 0.
\label{eq:hermite}
\end{equation}
Euler's and trapezoidal integration rules are computationally less expensive
than the one given in (\ref{eq:hermite}).
Therefore, the integration scheme influences the accuracy of the
solution and also the possibility of the solver of finding a feasible
solution.
To mitigate these issues a more refined discretization, i.e. more nodes, might
be necessary. In return however, the size of the problem and
potentially the time required to find a solution is increased.
Depending on the ability of the solver to handle the sparsity in the structure
of the problem, this increase in complexity may be partially absorbed by the
solver itself.
Still, there is a compromise between simplicity of the integration rule, number
of nodes and efficiency of the solver to be made.
%

\subsection{NLP solvers}

Nonlinear programming solvers are in a mature state of development given the
vast experience of the field of mathematical optimization.
%
%
Here we evaluate two representative instances of NLP solvers. The first solver
in our comparison is SNOPT which is based on Sequential Quadratic Programming.
The second solver is IPOPT using the Interior Point Method. The goal of the
analysis is to evaluate the performance of both approaches within a robotics
application.
The complete description of the solvers and the algorithms lies outside of the
scope of this paper. This section presents some features aiming to describe the
main differences between them and to identify the nature of any difference in
their performances.
A complete comparison between these classes of solvers can be found in
\cite{Betts2002}.

\subsubsection{SNOPT (Sparse Nonlinear OPTimizer)}
The basic structure of this implementation of the SQP algorithm involves major
and minor iterations.
Major iterations advance along a sequence of points ${y_h}$ that satisfy the set
of linear constraints in $\zeta(y_h)$ and $g(y_h)$.
These iterations converge to a point that satisfies the remaining nonlinear
constraints and the first-order conditions of optimality \cite{Gill2005}.
The direction towards which the major iterations move is produced by solving a
QP subproblem.
Solving this subproblem is an iterative procedure by itself (i.e. the minor
iterations), based on a Newton-type minimization approach.

An important characteristic of all SQP algorithms is that they are `active set'.
Roughly speaking, this means that during the iterative procedure all the
inequality constraints play a very explicit role as the QP subproblem must
estimate the active set in order to find the search direction. 


\subsubsection{IPOPT} 
This algorithm also depends on a Newton type subproblem. Nevertheless,
inequalities are handled in a different manner. A barrier function is used to
keep the search as far as possible from the bounds of the feasible set. The
barrier parameters change along iterations, allowing proximity to the adequate
constraint.    

%
%
%
%

\section{Evaluation Procedure}
\label{sec:problem_set}

Considering the elements presented in the previous section, there are several
decisions to make in order to put all the pieces together to solve a robot
motion planning problem using direct transcription.

What integration method to use? What is the influence of the number of nodes?
Given a complex robotic problem (i.e. high state dimensionality, nonlinear
dynamics and path constraints), how long does it take to find a solution? Can
this technique be implemented online? What type of solver is more favorable for
TO in robotics? Can the resulting trajectory be implemented on a real robot?

We address these questions measuring the performance of different configurations
of the method against a robot motion benchmark.

\subsection{Benchmark robot and tasks}
\label{sec:benchmark}
The robot used for the benchmark is the ball balancing robot Rezero
\cite{Fankhauser2010}, a so called ballbot. These kind of robots are essentially
3D inverted pendulums and hence are statically unstable, underactuated and
non-minimal phase systems.
%
Due to this instability there is a complex interaction between the requirements
for stable control and satisfying constraints.

\begin{figure}
	\centering
	\includegraphics[width=0.95\columnwidth,
keepaspectratio=true]{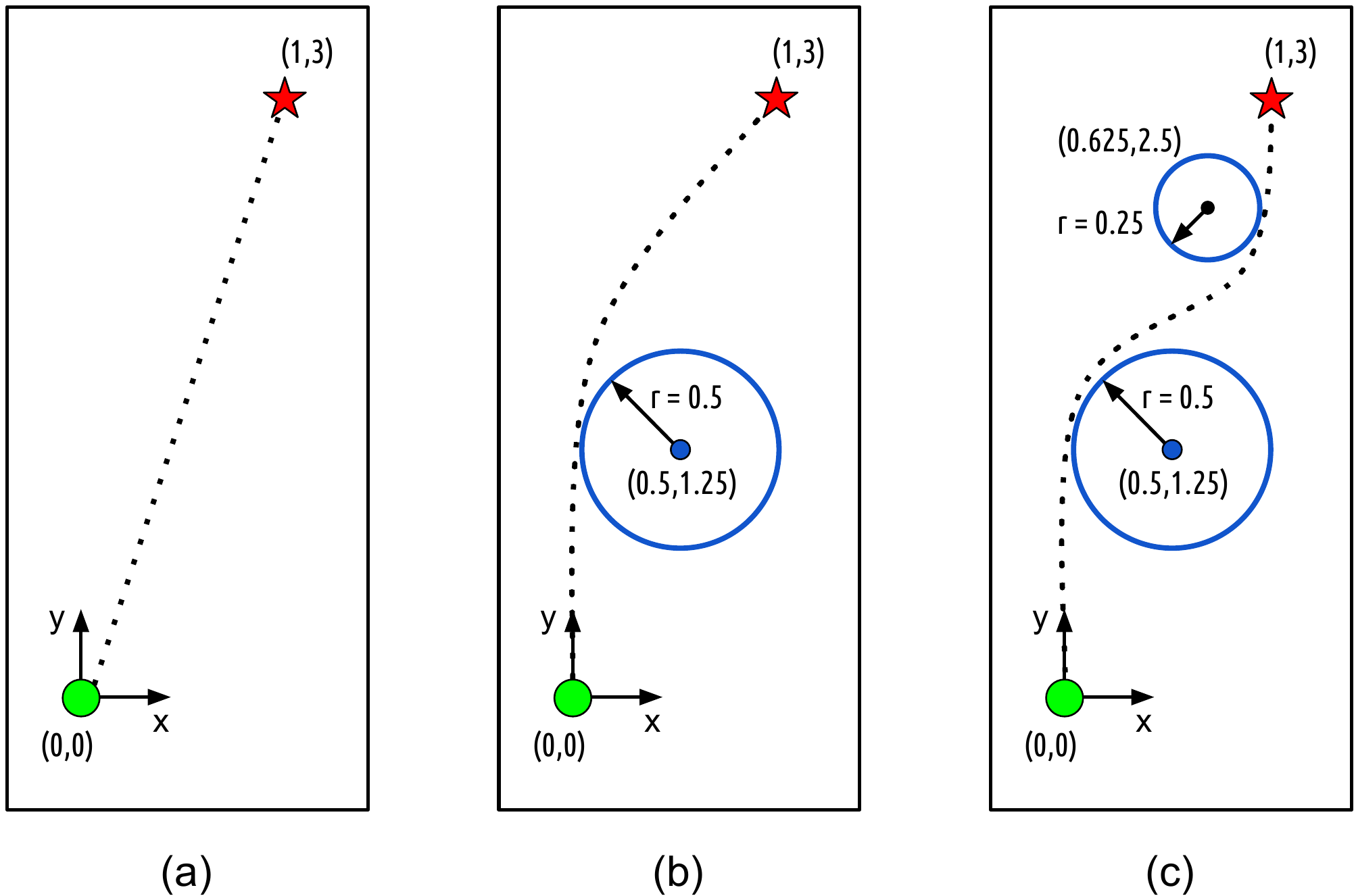}
	\caption{Illustration of the benchmark tasks (top view). Starting from the
origin the robot should reach the goal (star mark). Circle obstacles (blue) are
included to increase the problem complexity.  (a) Go-to task, (b) one-obstacle
task and (c) two-obstacles task. Position and size of the obstacles are
expressed in meters.}
\label{fig:benchmarktask}
\end{figure}

The benchmark consists of a set of three different motion tasks.
Each task is a variation of a go--to task, i.e. the robot is supposed to move
from an initial to a final spatial location avoiding fixed obstacles (see
Fig.~\ref{fig:benchmarktask}). State and control trajectories are not
pre-specified and should result from the minimization of a cost function as well
as from the satisfaction of the boundary and path constraints. Details on the
implementation of the cost function and constraints are presented in Section
\ref{sec:robot_model_and_control_scheme}.
%

\subsection{Problem components}
The benchmark compares the performance of different configurations of the
direct transcription method. Table \ref{tb:schemevariables} summarizes the
variables evaluated in this study. 

Regarding the complexity of the task, each path constraint added to
(\ref{eq:constraint_3}) has to be verified at each node in (\ref{eq:NLP}),
augmenting the size and density of the gradients required to find an optimal
solution. This complexity is also increased when bounds on the state and control
are considered. This is the case for the hardware
experiments presented later in Section \ref{sec:hardware_experiments}. Here we
measure the effect of including such bounds to the two-obstacle task, in a
problem complexity labeled as `two-obstacle-bounds'.

At the same time, the performance also depends on the initial values of the
decision variables of the NLP. Here we evaluate three different initialization
methods, namely Zero, Linear and Incremental. In the first method all variables
are naively set to zero before the optimization begins. The linear method
initializes only the variables corresponding to the states of the
robot, using a linear interpolation between the initial and desired final
states. Finally, the incremental method uses the optimal solution of a simpler 
task 
(e.g. less constrained) as initial values for all variables. 

\begin{table}
\centering
\caption{Analyzed Components of a Direct Transcription Problem}
    \begin{tabular}{ | l | l |}
    \hline
    Integration Scheme & Trapezoidal \\
    					& Hermite-Simpson \\ \hline
    Problem Complexity  & Go-to \\ 
    					& One-obstacle\\
    					& Two-obstacle\\
    					& Two-obstacle-bounds\\
    					\hline
    NLP Initialization	& Zero \\
    					& Linear interpolation \\
    					& Incremental \\ \hline
    Number of Nodes		& \\ \hline
	NLP Solvers			& SNOPT \\
						& IPOPT \\ 
    \hline
    \end{tabular}
    \label{tb:schemevariables}
\end{table}

\subsection{Performance criteria}
\label{sec:performance_criteria}
The ultimate goal of TO is to obtain state and control trajectories that can be
implemented on a real robot, such that a certain task is accomplished while
an objective function is minimized. Accordingly, we measure the performance of the
methods in terms of accuracy, solving time and quality of the solution. 

\paragraph{Accuracy} The system dynamics are approximated using a piecewise
polynomial function. A solution is only valid if such an approximation is
accurate. In case of perfect accuracy, almost no stabilizing feedback control would be
required during simulation.
We measure the accuracy of the solution using the mean squared error between the
planned and the actual total control trajectories measured when the complete
system is simulated (i.e. an indication of the divergence between planned
$(u^*)$ and total $(u)$ control signals as defined in Section
\ref{sec:robot_model_and_control_scheme} and represented in Fig.
\ref{fig:control_scheme}).

\paragraph{Running time} Ideally, the motion planner should be fast enough to
compute trajectories online.
However, depending on the configuration, different number of
iterations might be required to solve the NLP. Here we use the \textsl{total
time} required to find a feasible solution as indicator of performance. In
contrast to the number of iterations, it provides an absolute value that can be
used to compare among solutions obtained with different solvers and integration
schemes.

\paragraph{Quality of the solution} The conditions of optimality for the
continuous trajectory optimization problem are given by the Pontryagin's minimum
principle. In direct transcription such conditions are approximated by the
Karush-Kuhn-Tucker (KKT) necessary conditions for the corresponding NLP (see 
e.g.,
\cite{Betts1998}). Continuous and discrete conditions converge as $N \rightarrow
\infty$ and $\Delta T \rightarrow 0$. In this paper the KKT conditions are
satisfied by any given solution with a tolerance of $10^{-6}$. We compare the
\textsl{quality} of the solutions in terms of optimality using the resulting
value of the objective function.

\section{Robot Model and Control Scheme}
\label{sec:robot_model_and_control_scheme}

\subsection{Robot model}
The non-linear model of the ballbot dynamics is described in
\cite{Fankhauser2010}.
The robot is modeled as two rigid bodies, the torso of the robot and the ball.
The two bodies are linked by three actuators. In our model, we neglect wheel
dynamics and we assume that no slip or friction losses occur.
The state vector of this ballbot is defined as
\begin{equation}
x = [\phi_{r} \quad \dot{\phi}_{r} \quad \phi_{p} \quad \dot{\phi}_{p} \quad
\phi_{y} \quad \dot{\phi}_{y} \quad \theta_x \quad \dot{\theta}_x \quad \theta_y
\quad \dot{\theta}_y]^T,
\end{equation}
where $\phi$ and $\dot{\phi}$ represent the torso angles and velocities
in the roll, pitch and yaw directions $(rpy)$. 
Furthermore, the state includes the rotational angles of the ball $(\theta_x,
\theta_y)$ as well as their derivatives $(\dot{\theta}_x,\dot{\theta}_y)$,
representing the ball position with respect to the initial state. The control
actions are defined by the wheels' input torques $u = [\tau_1 \, \tau_2 \,
\tau_3]^T$.

The desired goal state for all tasks is given by,
\begin{equation*}
x_g = \left[ 0 \quad 0 \quad 0 \quad 0 \quad 0 \quad 0 \quad 8 \quad 0 \quad 24
\quad 0  \right]^T
\end{equation*}
where the goal angles of the ball ($\theta_x = 8$ rad, $\theta_y = 24
$ rad) correspond to the desired spatial location ($x = 1$ m , $y = 3$ m),
given that the radius of the ball is $r=0.125$ m.

\subsection{Cost function, terminal conditions and total trajectory time
constraint}
%
%
A quadratic function is used for the intermediate cost, i.e.,
\begin{equation}
\psi(x_k,u_k)  =  x_k^TQx_k + u_k^T R u_k.
\end{equation}
Where $\bar{x}_N = x_N - x_g$ is the difference between the last state of the trajectory $x_N$ and the desired goal state $x_g$.
The diagonal cost matrices $Q,R$ weigh the contribution of the state and
control to the total cost.
%
%
Setting those values of $Q$ to zero that correspond to the rotational angles of 
the ball, tends to reduce the amount of time
the robot is not in the upright pose, and at the same time use as little
torque as possible along the path.

In direct transcription initial conditions are assumed to be fixed.
Additionally, in the simulations and experiments presented in this paper, the
desired goal state is added as a hard terminal condition to the NLP, (i.e., 
$x_N =
x_{g}$). In case this condition is not satisfied the problem is
considered infeasible. It is important to note that by using this hard constraint
the use of a final cost function $\Psi(x_N)$ can be avoided.

Finally, the time between nodes $(\Delta T)$ is chosen to be constant and included as a decision
variable. A constraint is added bounding the total trajectory time $(T)$, 

\begin{equation}
T_{min} \leq (N-1) \cdot  \Delta T \leq T_{max},
\end{equation}
where the bounds are given as a task specification. For this study we define that
motions should be completed between $T_{min} = 1.0$ s, and $T_{max} = 3.5$ s. 

\subsection{Trajectory stabilization}
\label{subsec:trajectory_stabilization}
Despite using an elaborated integration scheme or a very small discretization
time, solutions obtained using direct transcription are based on an
approximation of the system dynamics. 
Even in simulation very small numerical integration errors push the system
towards instability.
Therefore, it is necessary to stabilize the system.

\begin{figure}
\centering
\includegraphics[width=0.80\columnwidth,keepaspectratio=true]{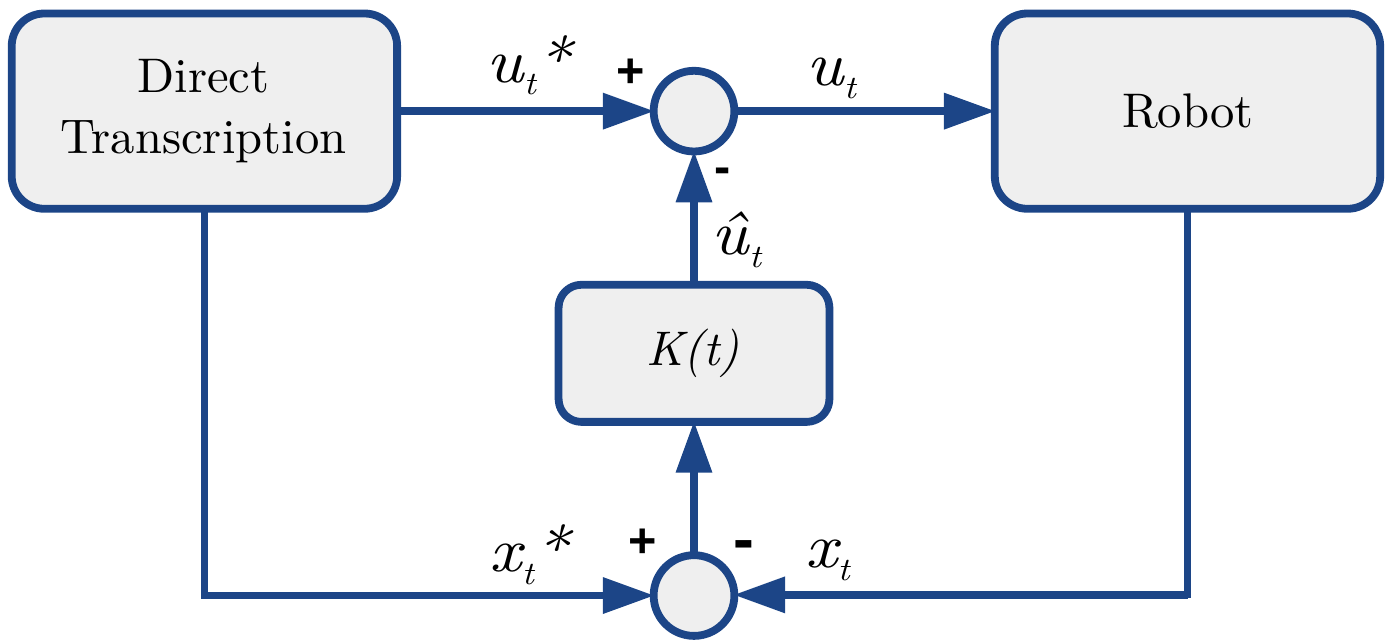}
\caption{Complete control scheme. Feed forward control action is provided by the
TO planner and a feedback stabilization control is computed using a time varying
gain matrix, $K(t)$, obtained via TVLQR.}
\label{fig:control_scheme}
\end{figure}

In this work we use Time Variant Linear Quadratic Regulator - TVLQR (see e.g.,
\cite{Tedrake2014}) to obtain a feedback controller. The derivation of this
optimal controller is out of the scope of this paper. However, a descriptive
explanation is provided for completeness: The dynamics of the system are
linearized around the continuous optimal trajectory $\left(x^*_t,u^*_t\right)$
using a Taylor approximation of (\ref{eq:systemdynamics}). The state and control
trajectory errors are defined as $\hat{x}_t = x^*_t - x_t$ and $\hat{u}_t =
u^*_t - u_t$.
Thus, its dynamics are described as a linear time varying system
\begin{equation}
\dot{\hat{x}}_t = A(t)\hat{x}_t + B(t)\hat{u}_t.
\label{eq:tvsystem}
\end{equation}
 The state and input matrices are given by
\begin{eqnarray*}
A(t) & = & \left.\frac{\partial f\left(x,u\right)}{\partial x}
\right\|_{x^*_t,u^*_t}\\
B(t) & = & \left. \frac{\partial f\left(x,u\right)}{\partial
u}\right\|_{x^*_t,u^*_t}.
\end{eqnarray*}
By solving the TVLQR problem for the system in
(\ref{eq:tvsystem}), a matrix of time dependent optimal feedback gains $K(t) \in
\mathbb{R}^{m\times n}$ is obtained. The resulting control scheme is illustrated
in Fig. \ref{fig:control_scheme}.

\section{Results}
\label{sec:results}

The results presented in this section were obtained using a standard laptop
computer with 2.4Ghz Intel Core i5 (dual core) processor with 4GB 133Mhz DDR3
RAM. Feasibility tolerance of the NLP solvers was set to $1\times10^{-6}$,
meaning that all the constraints, including the terminal condition $(x_N =
x_{g})$, are satisfied with an error bounded by this value.

\subsection{Optimal solution and stabilizer}

Here we present the resulting state and control trajectories for the goto task
obtained with a typical configuration: 100
Nodes, Hermite-Simpson integration scheme, linear initialization and SNOPT
as NLP solver.

\begin{figure}
\centering
\includegraphics[width=0.85\columnwidth,keepaspectratio=true]{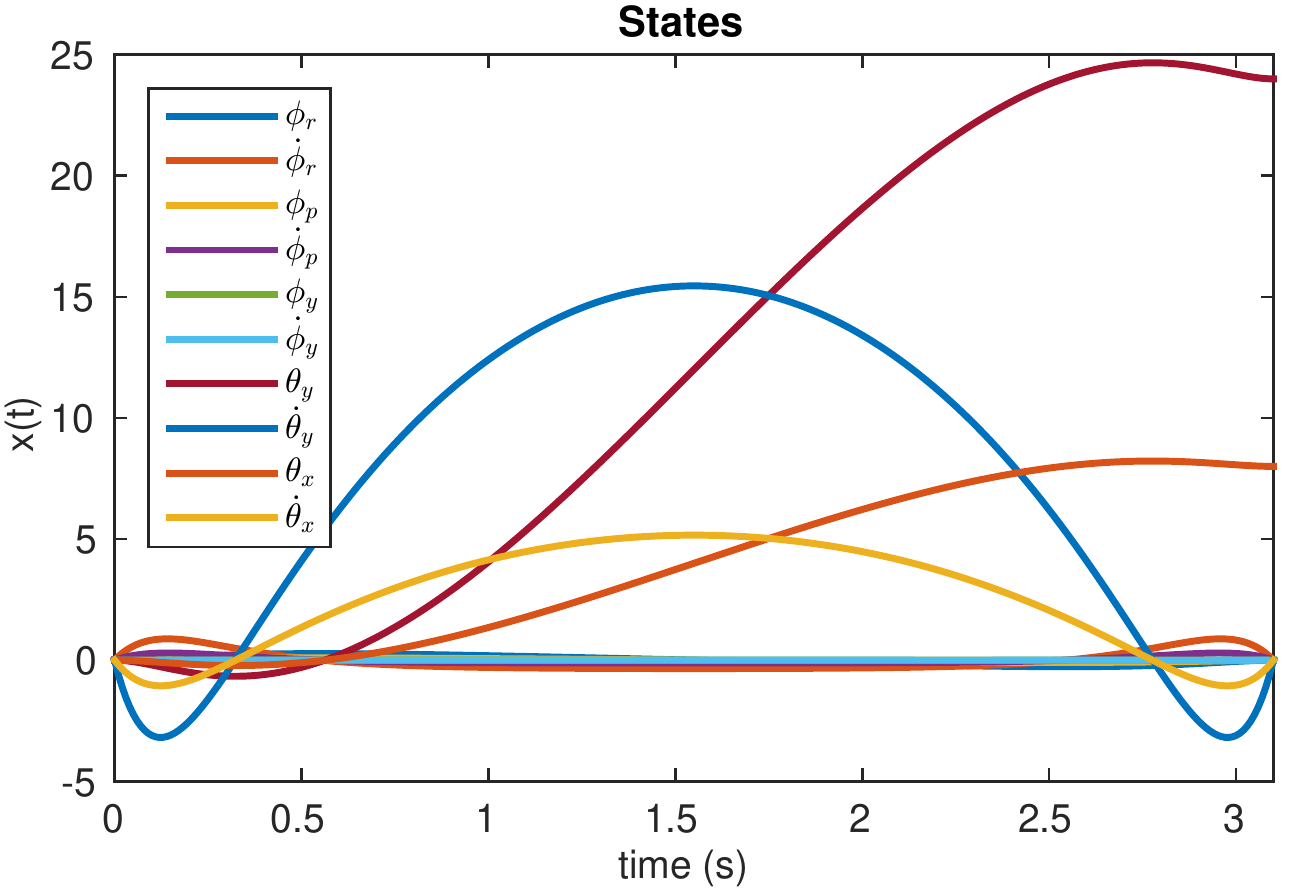}
\caption{State trajectories for the goto task. Solution found using $N=100$
nodes, linear initialization, Hermite-Simpson integration and SNOPT solver. The
optimal solution also includes the time increment between nodes $\Delta T =
0.03131$ s. All velocities and body angles converge to zero.}
\label{fig:x_solution}
\end{figure}

\begin{figure}
\centering
\includegraphics[width=0.85\columnwidth,keepaspectratio=true]{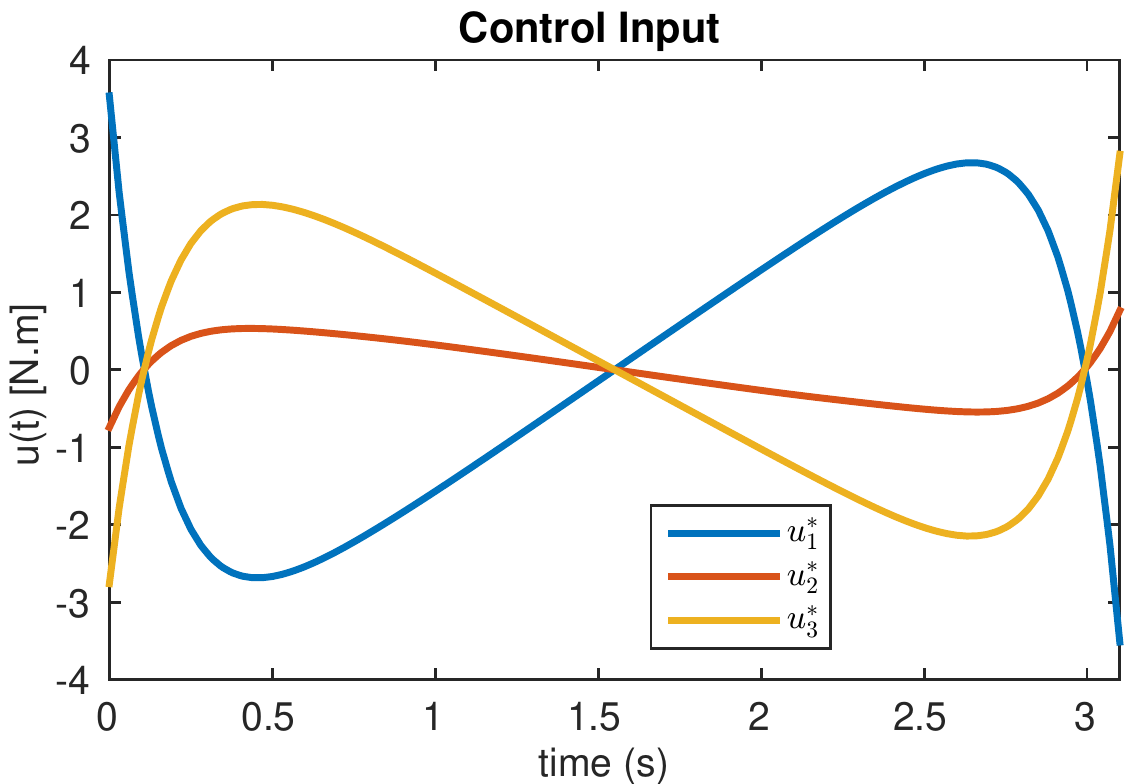}
\caption{Control trajectories for the goto task. Solution found using $N=100$
nodes, linear initialization, Hermite-Simpson integration and SNOPT solver.}
\label{fig:u_solution}
\end{figure}

\begin{figure}
\centering
\includegraphics[width=0.85\columnwidth,keepaspectratio=true]{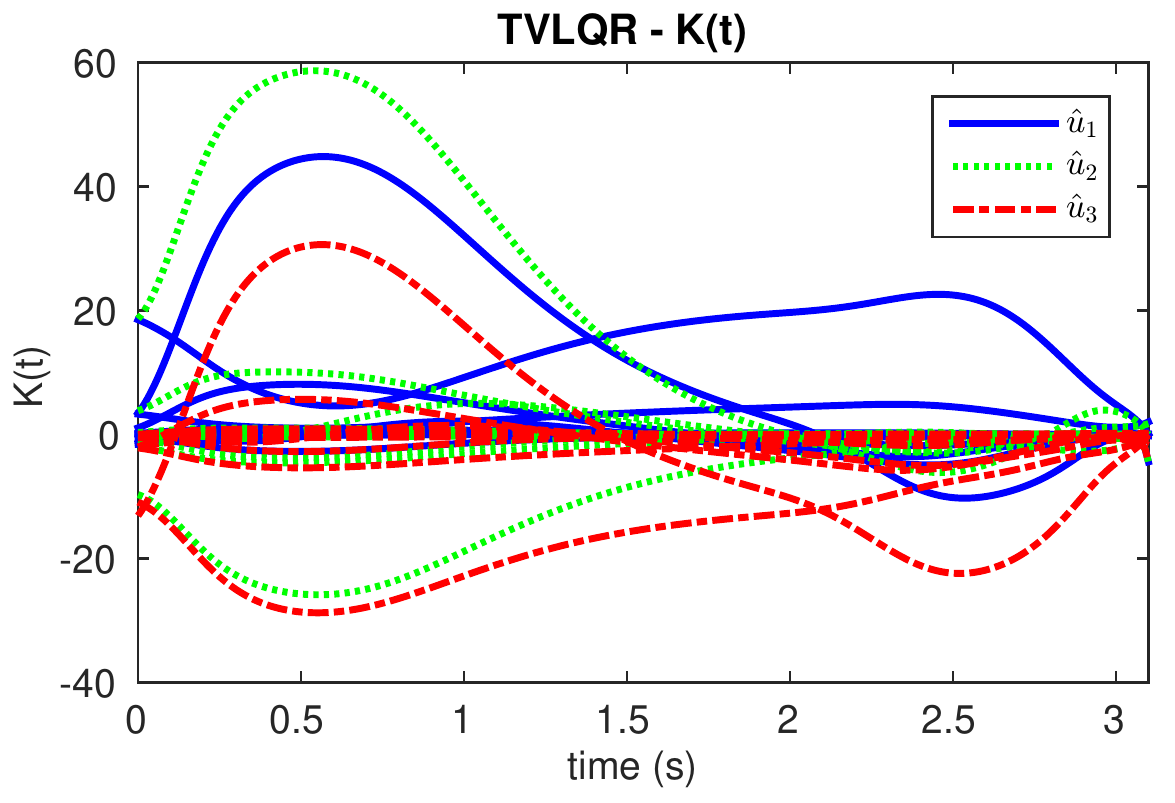}
\caption{Feedback gains for stabilizing the goto task. Gains are grouped
(line/color) by the corresponding control. Each feedback control $(\hat{u}_i)$
is computed adding the weighted contribution of the state error ($\hat{x}$).} 
\label{fig:K_solution}
\end{figure}

Fig. \ref{fig:x_solution} presents the optimal trajectories for the states. The desired state is reached in $T=3.1$ s. Similarly, Fig.
\ref{fig:u_solution} shows the control trajectories and Fig.
\ref{fig:K_solution} shows the corresponding stabilizer gains obtained using
TVLQR.

\subsection{Task complexity and initialization strategy}

In this part, the performance of direct transcription for tasks of diverse 
complexity
as well as the influence of different NLP initialization strategies are
evaluated.

\begin{figure}
\centering
\includegraphics[width=\columnwidth,keepaspectratio=true]{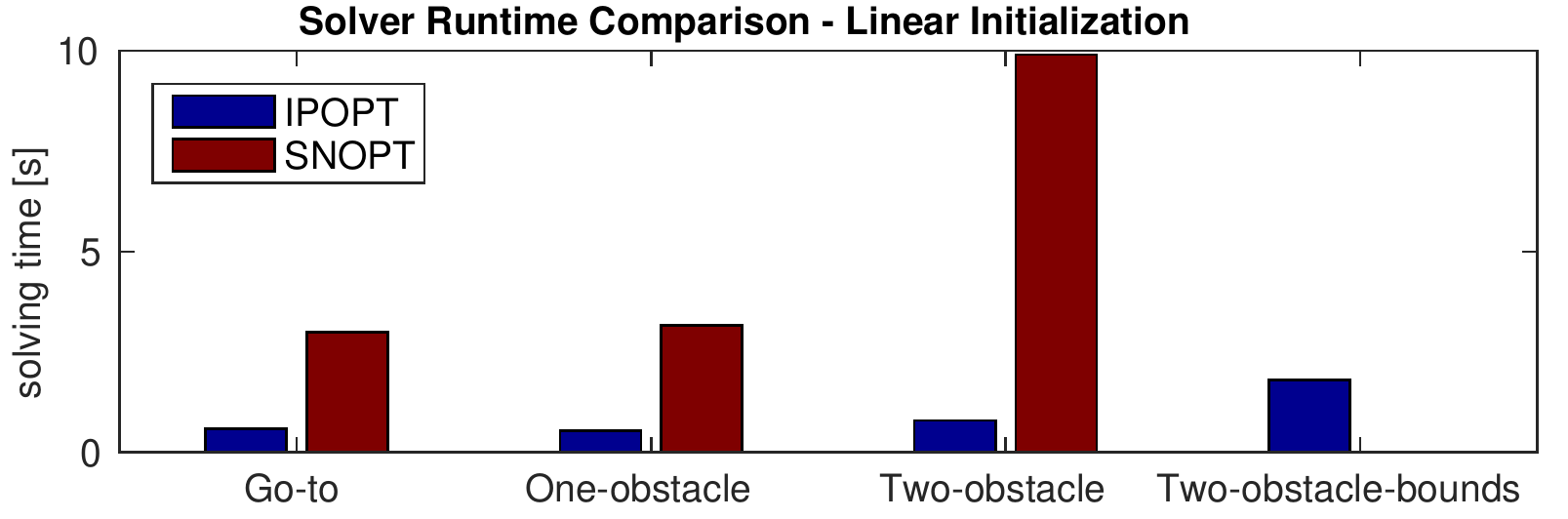}
\caption{This plot compares the solving time for the different tasks and solvers
using linear initialization method. Number of nodes $(N=100)$ and integration
scheme (trapezoidal) are constant. SNOPT is not able to find a solution for the
real robot task given those conditions.} 
\label{fig:time_lin_init}
\end{figure}

Firstly, we observed that the two-obstacles task cannot be solved using zero
initialization regardless of the number of nodes, integration scheme or solver
used. On the other hand, such a task can be solved using linear initialization.
Fig. \ref{fig:time_lin_init} shows the solving time required by the solvers to
find a solution to the different tasks using the same scheme (100 nodes,
trapezoidal interpolation and linear initialization).
It can be seen that the complexity of the task affects the time required by the
solvers. This is more significant in the case of SNOPT.
Moreover, for the two-obstacle-bounds complexity, SNOPT does not find a feasible
solution before reaching the maximum number of iterations ($1000$),
irrespective of the integration scheme or number of nodes used.

\begin{figure}
\centering
\includegraphics[width=\columnwidth,keepaspectratio=true]{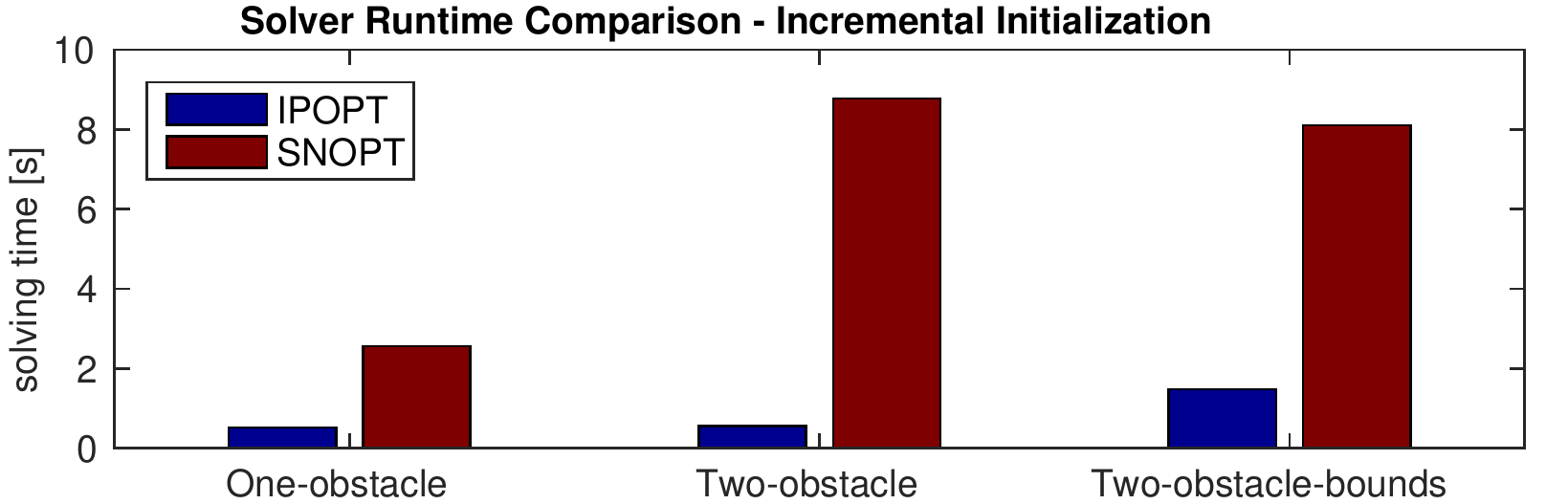}
\caption{This plot compares the solving time for the different tasks and solvers
using incremental initialization. Number of nodes $(N=100)$ and integration
scheme (trapezoidal) are constant.} 
\label{fig:time_homo_init}
\end{figure}

In Fig. \ref{fig:time_homo_init} solving time for the different tasks and
solvers are shown for the case of incremental initialization. Decision variables
are initialized with the solution of the previous task. As expected in
this type of initialization, the solving time is lower than those shown in
Fig.\ref{fig:time_lin_init}. Moreover, in this case SNOPT is able to find a
solution for the two-obstacle-bounds complexity. Apart from the variation on running
time and feasibility reported in this section, no significant differences were
observed on the accuracy or quality of the solution given different
initialization methods.

\subsection{Integration scheme}

In this section we show the influence of the integration scheme and solvers in
terms of run time, accuracy and quality of the solution. In order to obtain the
corresponding data we use the one-obstacle task using zero initialization.

\begin{figure}
\centering
\includegraphics[width=0.85\columnwidth,keepaspectratio=true]{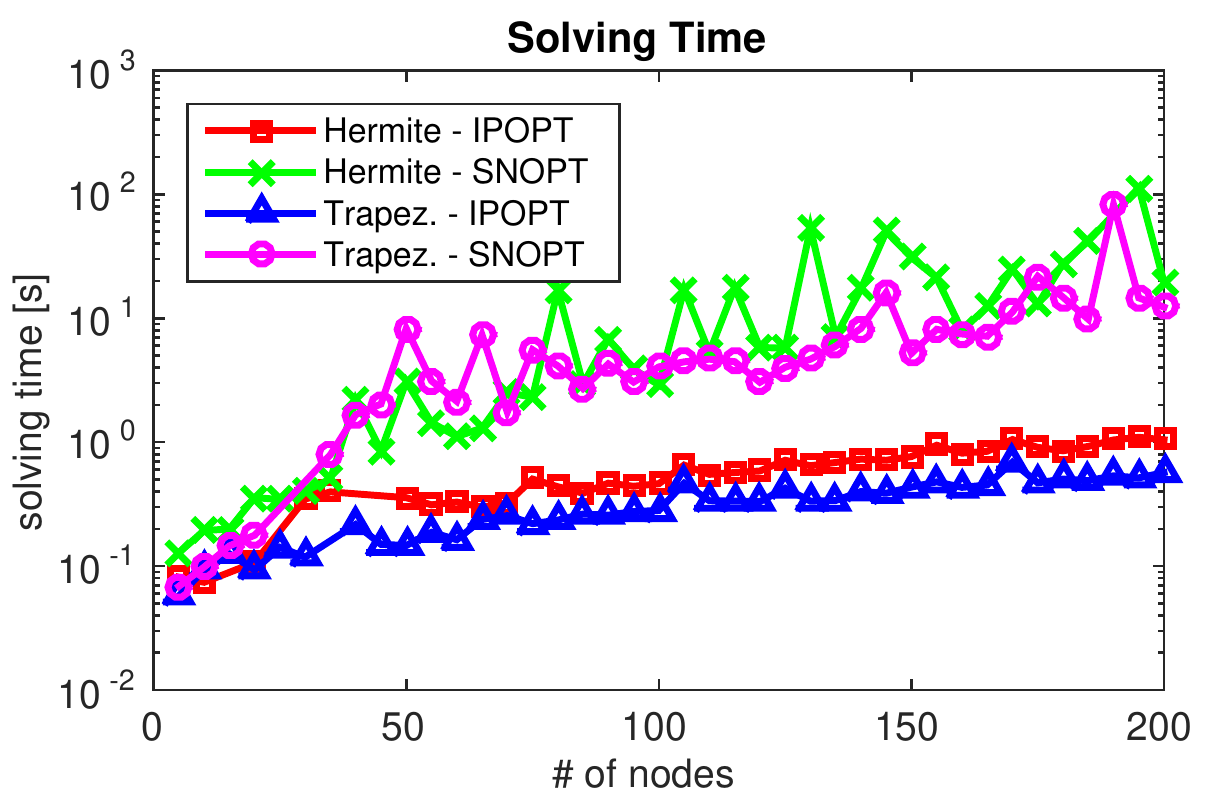}
\caption{This plot shows the runtime of the different combinations of
integration methods and solvers. For a low number of nodes, differences are
small. However, with an increasing number of nodes, the runtimes of SNOPT are
above those of IPOPT and also rise faster.}
\label{fig:cpu_time}
\end{figure}

Fig. \ref{fig:cpu_time} compares the runtime for the integration schemes in
combination with both solvers. While the absolute values of these measures are
highly dependent on the CPU capacity, this section is concerned with their
relative values. With a low number of nodes (i.e. $N < 20$) both solvers show
similar performance. However, with an increasing number of nodes the
computational complexity of SNOPT is higher and increases faster than the one
obtained with IPOPT. For both solvers, Hermite integration requires more CPU
time than trapezoidal integration, but the difference is not
significant.

\begin{figure}
\centering
\includegraphics[width=0.85\columnwidth,keepaspectratio=true]{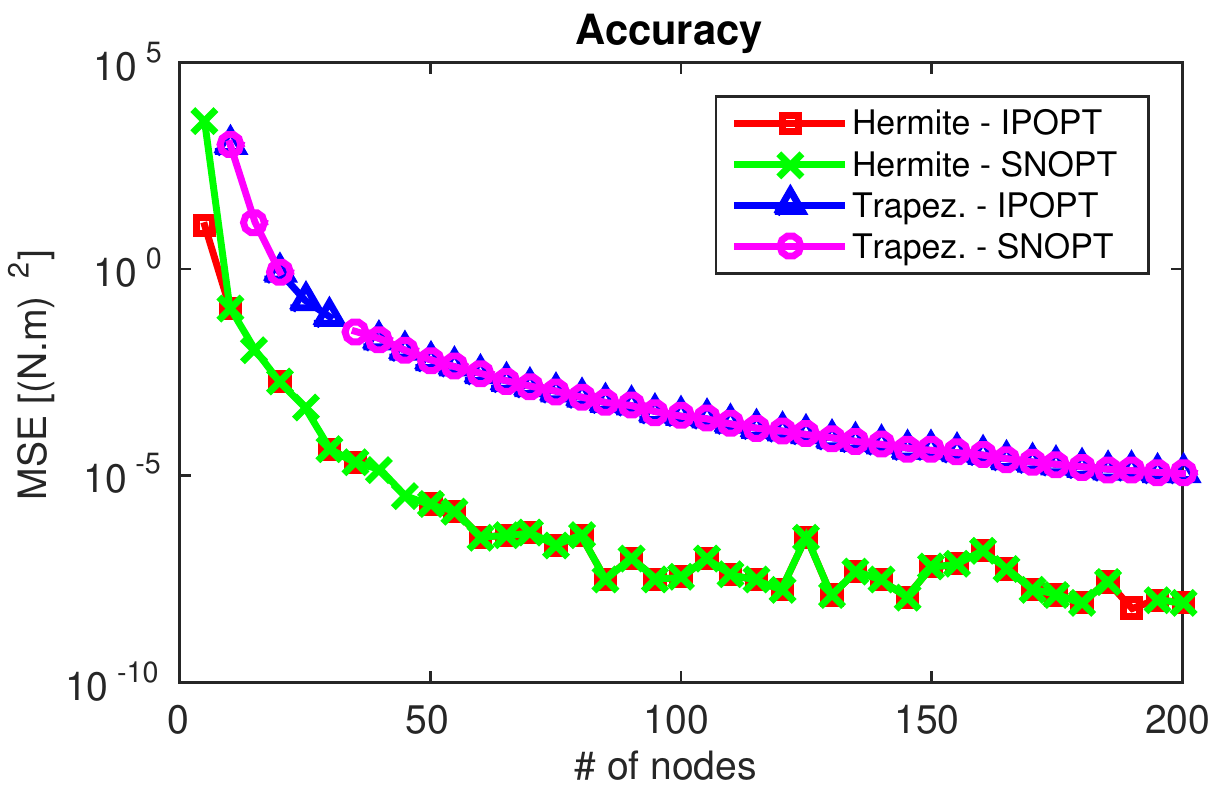}
\caption{This plot shows the accuracy of the solution for both solvers and
integration schemes. As indicator of accuracy we use the mean squared error
(MSE) between the planned and the total control required to simulate the task.
It can be seen that Hermite integration is significantly more accurate than
trapezoidal integration.}
\label{fig:integration_accuracy}
\end{figure}

Fig. \ref{fig:integration_accuracy} shows the divergence, in terms of mean
squared error, between the planned and total control trajectories (during
simulation). As it can be observed, the accuracy is independent of the solver.
However, Hermite integration produces an error which is by several magnitudes
lower than the one of trapezoidal integration.

\begin{figure}
\centering
\includegraphics[width=0.85\columnwidth,keepaspectratio=true]{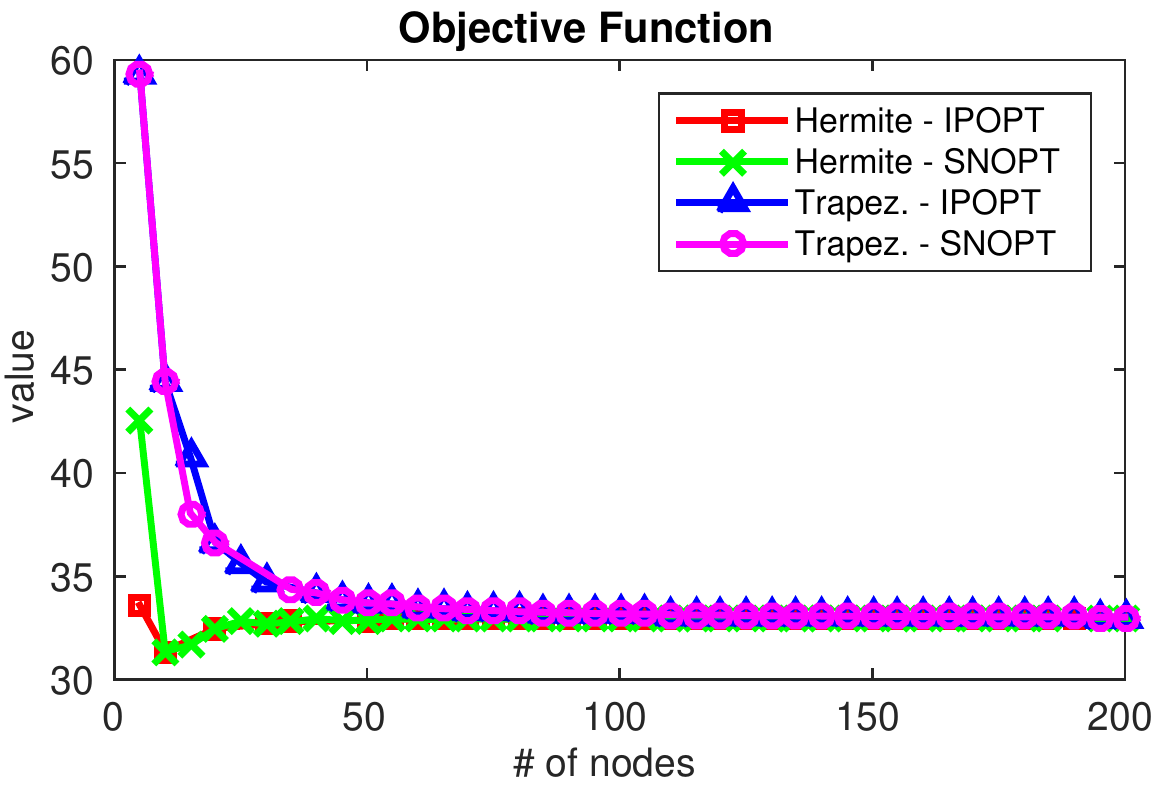}
\caption{This plot shows the objective function value for the different
combinations of integration methods and solvers. 
Especially for low number of nodes, Hermite integration clearly outperforms
trapezoidal integration. For more than 50 nodes, the differences vanish.}
\label{fig:objective_function}
\end{figure}

Fig. \ref{fig:objective_function} compares the value of the objective function
after optimization. This value can be understood as being inversely proportional
to the quality of the solution. The difference between the two solvers is
marginal. However, Hermite integration is able to achieve a significantly lower
objective function value, especially with a low number of nodes.

\subsection{Hardware experiments}
\label{sec:hardware_experiments}
\begin{figure}[tpb]
  \begin{minipage}[c]{0.3\columnwidth}
    \includegraphics[width=\columnwidth]{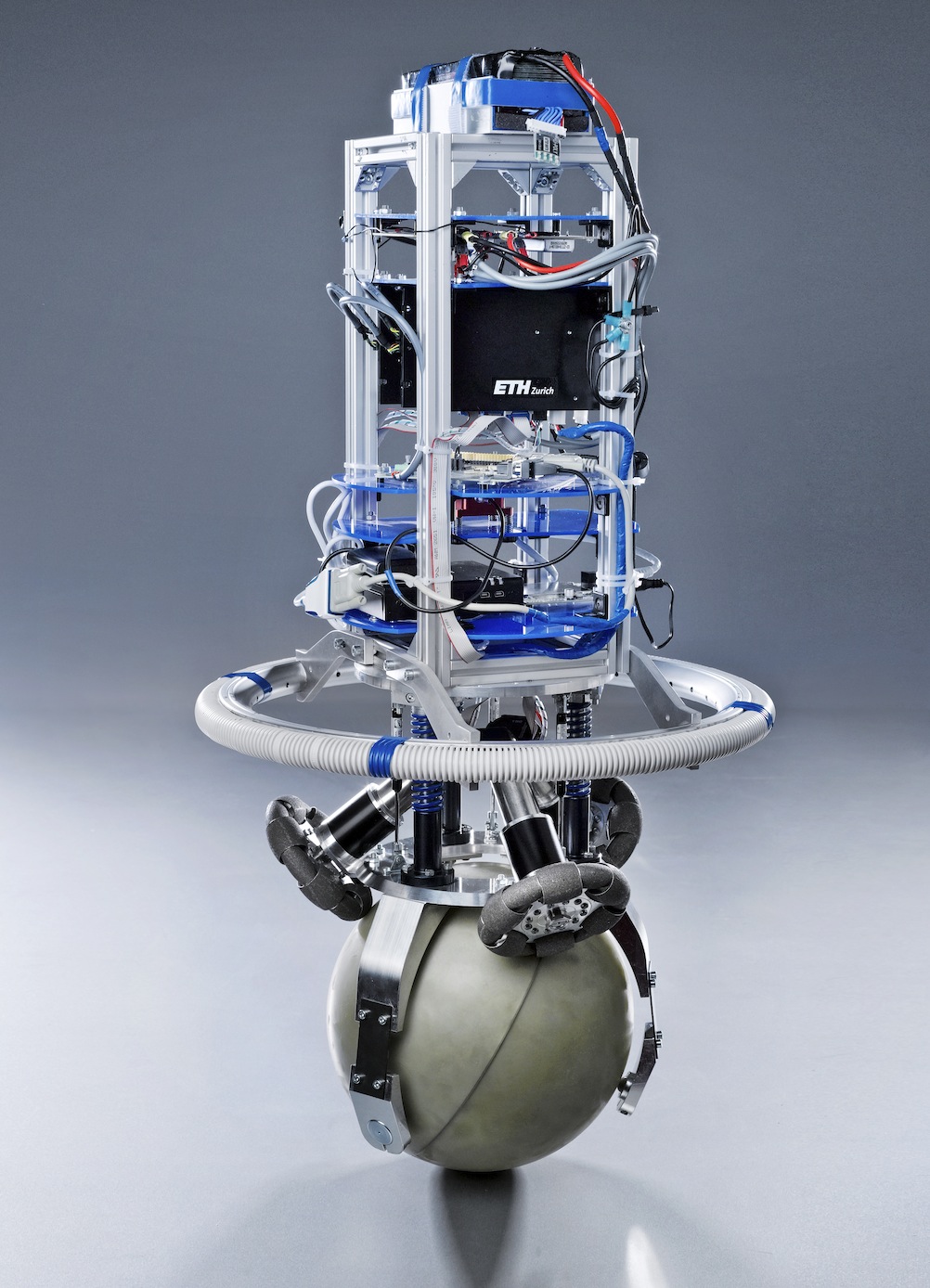}
  \end{minipage}\hfill
  \begin{minipage}[c]{0.60\columnwidth}
    \caption{
      Ball balancing robot. This robot consists of three major elements: the
ball, the propulsion unit and the upper body (torso). The ball is driven by
three omniwheels mounted on brushless motors with planetary gear heads. Through
optical encoders on the motors, the ball rotation is measured providing onboard
odometry for the robot.} \label{fig:rezero}
  \end{minipage}
\end{figure}

All the benchmark tasks shown in Fig. \ref{fig:benchmarktask} (including those
with path constraints) were applied to the real hardware in order to verify that
these approaches also hold on a physical system. Control inputs and states were 
bounded to more conservative values avoiding aggressive behaviors and to operate in a safer regime. Moreover, upper bound for the total trajectory time was relaxed ($T_{max} = 7$ s).
Solutions were found using 100 Nodes,  
Hermite approximation and IPOPT as a solver. The robot used in the
experiment is shown in Fig. \ref{fig:rezero} and a video with the experimental results can be found at
\url{https://www.youtube.com/watch?v=VGIROnFWgMw}. In the video it can be seen
that the robot executes the task synchronized with the simulation, following the
planned trajectories.


%

\begin{figure}[htbp]
\centering
\includegraphics[width=0.9\columnwidth,keepaspectratio=true]{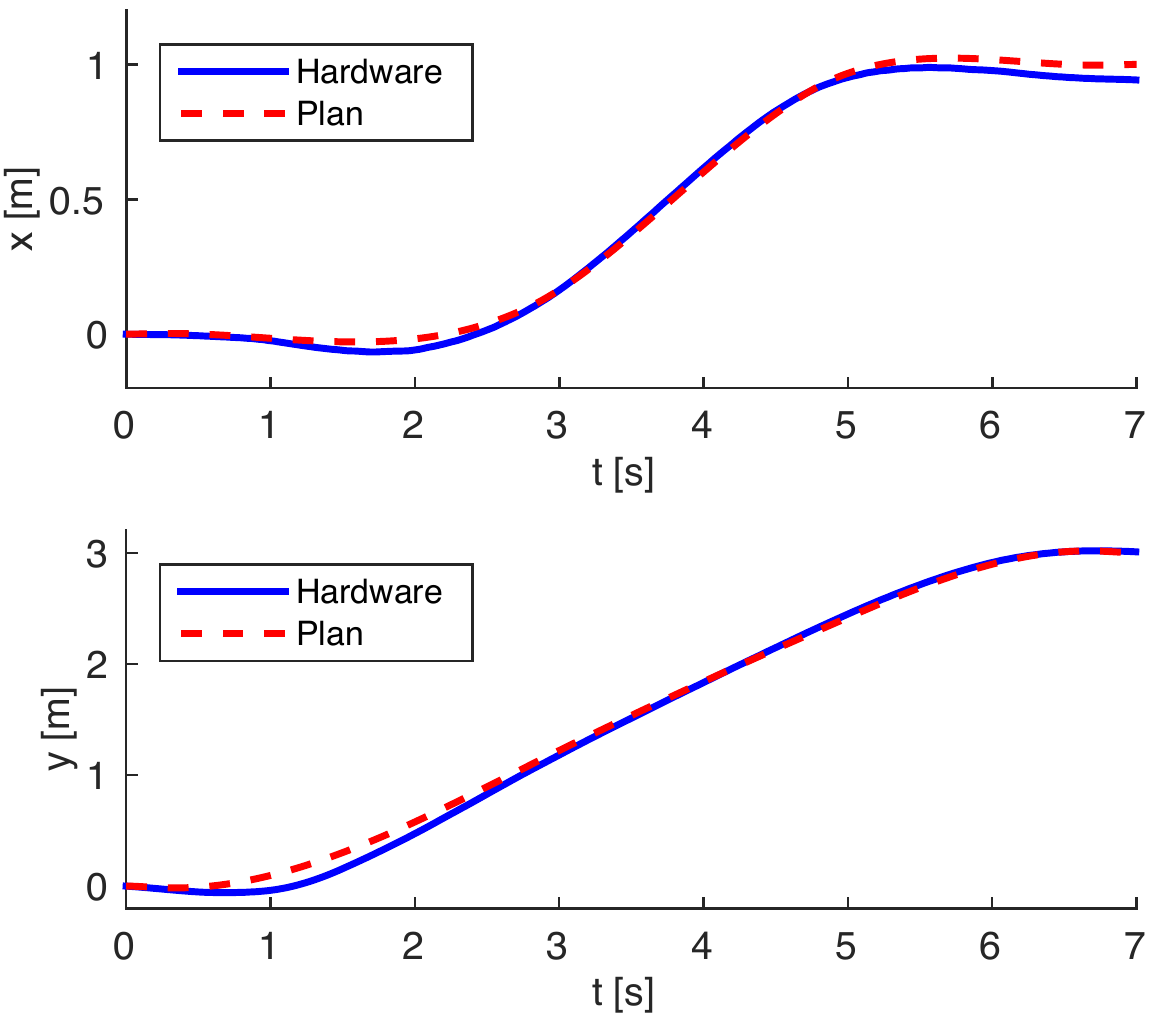}
\caption{Position of the robot in $x$ and $y$ direction during execution of the
two obstacles task. The blue solid line shows the measured
data while the red dashed lie shows the planned trajectory. Solution time for this configuration is $1.79$ s.}
\label{plot:slalom_xy}
\end{figure}

The $xy$ position of the robot during execution of all tasks matches the planned trajectory fairly well
with almost no deviation. Fig. \ref{plot:slalom_xy} shows the results for the
third task. Even though this is the most challenging task, it can be seen that
the tracking behavior is very good. In summary, the hardware tests show good
tracking performance, verifying that the optimized trajectory is dynamically
feasible and that the tracking control is able to reject disturbances and model
inaccuracies very well.

\subsection{Validation on a different robot}

\begin{figure}[t]
\begin{minipage}[c]{0.6\columnwidth}
\centering
	\includegraphics[width=\columnwidth,keepaspectratio=true]{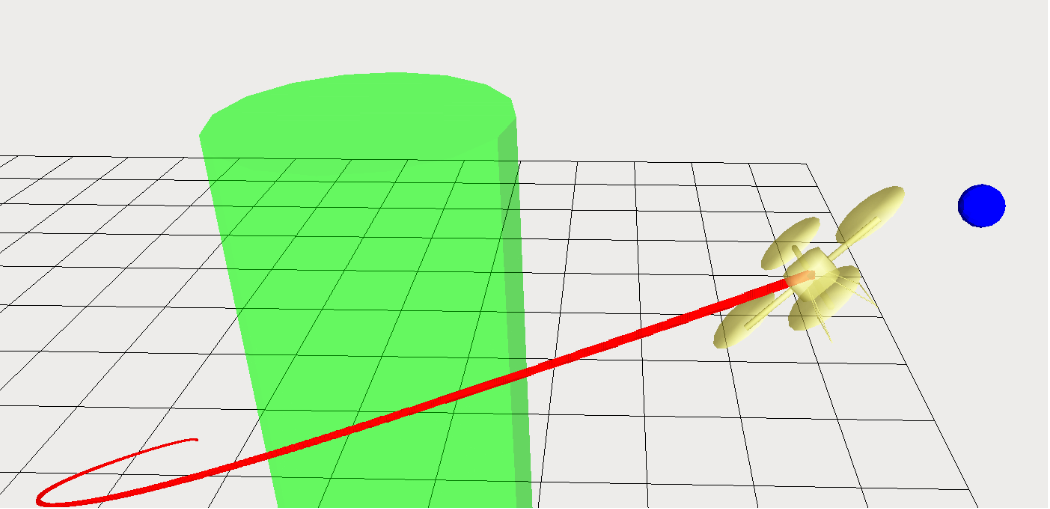}
	\end{minipage}\hfill
\begin{minipage}[c]{0.38\columnwidth}
  \caption{Quadrotor validation task. Blue dot represents the goal spatially located at $p = [5,5,5]$. A cylindrical obstacle ($r =1$) is centered at $x=2.5 , y=2.5$.}\label{image:quadrotor}
\end{minipage}
\end{figure}

Results presented above are validated by applying the same
methods on a simulated 3D quadrotor ($x \in \mathbb{R}^{12}, u \in
\mathbb{R}^{4}$) executing an obstacle avoidance task.
Fig. \ref{image:quadrotor} shows the trajectory followed by the quadrotor, starting from the origin and reaching the goal after avoiding
a cylindrical obstacle. The complete motion can be seen in the supplemental
video.

We observed that the influences of the type of solver, dynamic
constraints, initialization method and number of nodes are similar to the ones
shown on the ballbot experiments. For instance, Table \ref{tb:verification}
shows 
the time required to solve this task for different configurations using zero initialization strategy.
The same patterns can be observed: Solution time increases with the number of
nodes and IPOPT is faster than SNOPT given the same integration scheme. As
observed in Fig. \ref{fig:cpu_time}, there is no clear pattern regarding
solution times required by SNOPT when using different integration schemes. 
These results support the hypothesis that the conclusions drawn from the evaluation
presented in this paper are not robot specific but that they can be extended to
other platforms.

\begin{table}
\caption{Verification Results - Solution Time (\MakeLowercase{s})}
\label{tb:verification}
	\centering
	\begin{tabular}{|c|c|c|c|c|}
		\hline
		\multirow{2}{*}
		{N} &
			\multicolumn{2}{|c|}{IPOPT} &
			\multicolumn{2}{|c|}{SNOPT} \\
			\cline{2-5}
		 & Hermite & Trapez. & Hermite & Trapez. \\
		 \hline
		 50 & 0.785 & 0.240 & 18.53 & 0.423 \\
		 75 & 1.153 & 0.332 & 42.81 & 3.072 \\
		\hline
	\end{tabular}
\end{table}

\section{Analysis and Discussion}
\label{sec:discussion}
\setcounter{paragraph}{0}

\paragraph{NLP solver}
We observed that SNOPT requires more time in most of the comparisons.
This is even clearer for the case of the two obstacle task, where the difference
to IPOPT is significant.
This can be explained by the need of SQP methods in finding a solution within a
big active set of inequality constraints.
As pointed out in other comparison studies \cite{Betts2002}, IPM tend to be
faster in a NLP problem with many inequality constraints.
%
%
%
%
%
Moreover, given a good initialization, the IPM algorithm always outperformed the
SQP solver.
In terms of quality and accuracy of the solution both solvers acted similarly
under different schemes.

\paragraph{Integration scheme}
The integration scheme has less influence on the solution time than the solver.
However, selecting a Hermite-Simpson approximation has a considerable impact on
the accuracy of the solution. 
As expected, the error decreases with the number of nodes, but this is not
significant with respect to the better approximation reached with
Hermite-Simpson.
The quality of the solution is also affected by the integration scheme.
Trapezoidal methods requires at least $50$ nodes to reach the same quality as
the one obtained with Hermite integration using $20$ nodes. This fact, together
with the running time results, where Hermite and trapezoidal schemes perform
similarly when using IPOPT, suggests that Hermite integration scheme should be
preferred.

\paragraph{Initialization}
The aspect of initialization is closely related to the performance of the
solver. 
We observed that direct transcription in general is considerably sensitive to the
initialization strategy adopted at the NLP stage.
Regardless of the solver, applying an adequate initialization policy
clearly reduces the time to obtain a solution. This is a fundamental aspect for 
applying direct
transcription online.
 
%
%
%
%
%


\paragraph{Hardware experiments}
It has been shown that TO with a matching controller
can be deployed on real hardware. The good match between simulated and
real behavior can be observed both from the figures presented in this paper as
well as from the video attachment overlaying simulated and physical robot. Small
deviations between planned and real trajectories can be explained by model mismatches
and sensor noise.
These experiments validate the applicability of TO to real hardware
and they are one of very few examples of using this technique on real
underactuated robots. 

\paragraph{Solving time}
%
Ideally, the robot itself should be able to obtain plans online.
Since the optimized trajectory is usually complemented with a stabilizing controller
TO is not  bound to run at the same rate or
in hard real-time. 
%
However, the solving time of TO should be fast enough to react to large
scale disturbances or unforeseen changes in the robot's environment. 
The maximally acceptable solving time depends on the dynamics of the system 
and on the task. 
For the given example in Figure \ref{fig:time_homo_init} the solving time is
about 5 to 50 cycles of the inner stabilizing control loop.
These solving times lie within a reasonable magnitude for online planning.

\paragraph{Failures during optimal trajectory search}

\begin{table}
\centering
\caption{Missing data points - configurations}
\begin{tabular}{c|c|c}
\hline
Integration Scheme & Number of Nodes & Solver \\
\hline
Hermite & 190 & SNOPT \\ 
Hermite & 15, 25, 40, 45 & IPOPT\\
Trapez. & 25,30 & SNOPT\\
\hline
\end{tabular}
\label{tb:solvererror}
\end{table}

Occasionally the solver might not be able to find a solution, declaring
numerical difficulties or because the maximum number of iterations has been
reached. As observed in Fig.~\ref{fig:cpu_time} and 
Fig.~\ref{fig:objective_function}, data is missing for a few 
configurations of solver, integration scheme and number of nodes. Such
configurations are reported in Table \ref{tb:solvererror}. Further investigation
is required to detect the fundamental reasons for these failures, as no clear
pattern of the variables analyzed in this study is observed. 

Finally, apart from the aforementioned cases, in 
Fig.~\ref{fig:integration_accuracy} five other points are also missing. Such 
cases
correspond to solutions for which the TVLQR method is not able to find a
feedback controller and therefore the accuracy as defined 
in
Section \ref{sec:performance_criteria} cannot be determined. This is observed 
for some configurations
with less than 30 nodes for both integration schemes and solvers. Further
investigation regarding the relation between direct transcription and the
solution of a TVLQR problem is required in order to identify the fundamental
reason of such effect.

\section{Conclusions} 
\label{sec:conclusions}

%

%
The comparisons performed in this work show clear patterns on the performance of
direct transcription methods with respect to the different integration methods,
solvers, initialization and problem types.
%

%
%

We showed a successful implementation of a direct transcription based motion
planning and control approach for a real robot. For a problem of medium
complexity, such as the one evaluated here, the solution convergences
sufficiently fast so online planning using these methods comes in reach.

%


\bibliographystyle{IEEEtran}
\bibliography{references}
\end{document}